\documentclass[11pt]{article}

\usepackage[margin=1in]{geometry}
\usepackage[T1]{fontenc}
\usepackage[utf8]{inputenc}
\usepackage{lmodern}
\usepackage{amsmath,amssymb,amsthm,mathtools}
\usepackage{hyperref}

\hypersetup{
  colorlinks=true,
  linkcolor=blue,
  citecolor=blue,
  urlcolor=blue
}

\newtheorem{theorem}{Theorem}[section]
\newtheorem{conjecture}[theorem]{Conjecture}
\newtheorem{problem}[theorem]{Problem}
\newtheorem{proposition}[theorem]{Proposition}
\newtheorem{lemma}[theorem]{Lemma}
\newtheorem{corollary}[theorem]{Corollary}
\theoremstyle{remark}
\newtheorem{remark}[theorem]{Remark}

\DeclareMathOperator{\diam}{diam}
\DeclareMathOperator{\Inv}{Inv}
\DeclareMathOperator{\Cay}{Cay}
\DeclareMathOperator{\spanop}{span}
\newcommand{\F}{\mathbb{F}}
\newcommand{\K}{\mathbb{K}}

\title{Edge-Number Bounds for the Inversion Diameter of Graphs}
\author{Jiawen Bo\thanks{School of Mathematical Sciences and LPMC, Nankai University, Tianjin 300071, China. Emails:\{bojiawen,xinyan\}@mail.nankai.edu.cn},
Anqi Li\thanks{School of Mathematical Sciences, Nankai University, Tianjin 300071, China. Email: anqili@nankai.edu.cn.},
Xiaopan Lian\thanks{Center for Combinatorics and LPMC, Nankai University, Tianjin 300071, 
China. Email: Lian@nankai.edu.cn.}, 
Xin Yan\footnotemark[1]
}
\date{\today}

\begin{document}
\maketitle

\begin{abstract}
The inversion of a set $X$ of vertices in an oriented graph reverses every arc with both endpoints in $X$.
The inversion graph $I(G)$ of a graph $G$ has the labelled orientations of $G$ as its vertices, two orientations being adjacent when a single inversion transforms one into the other, and the inversion diameter $\diam(I(G))$ is its diameter.
Answering a question of Havet, H\"orsch and Rambaud, we prove the  bound in terms of edge number $\diam(I(G)) \le 2\sqrt{|E(G)|}$,
and we complement it with a lower bound $\diam(I(G)) \ge \frac{|E(G)|}{|V(G)|}$
obtained by viewing $I(G)$ as a Cayley graph on $\F_2^{E(G)}$. We further refine the upper bound for bipartite graphs $G$ by showing   
 $ \diam(I(G))\le \max\left\{\rho,
  \left\lceil\log_2\bigl(2+\sigma(2^{\rho-1}-1)\bigr)\right\rceil\right\}$
where the two parts of $G$ have maximum degrees $\sigma$ and $\rho$, respectively.

%Furthermore, we refine the bounds on inversion diamter of triangle-free graphs and lexicographic products.

% We further  show that  the least $\alpha$ with
% $\diam(I(G))\le \alpha\sqrt{|E(G)|}$ for all $G$ lies in $[\sqrt2,2]$, where the lower bound $\sqrt2$ is   witnessed by complete graphs and $\alpha=\sqrt2$ under the maximum-degree conjecture $\diam(I(G))\le \Delta(G)$.  
% For triangle-free graphs, we prove  
% \[
%   \diam(I(G)) \le \sqrt{2|E(G)|}+\frac12\log_2(2|E(G)|).
% \]
% Turning to the bipartite setting, we give short, computer-free proofs that 
% \[
%   \diam(I(G))\le \max\left\{\rho,
%   \left\lceil\log_2\bigl(2+\sigma(2^{\rho-1}-1)\bigr)\right\rceil\right\}
% \]
% when the two parts have maximum degrees $\sigma$ and $\rho$.
% Finally, for lexicographic products with edgeless graphs we prove
% $\diam(I(F[\overline{K_t}]))\le \diam(I(F))\cdot t$ for $F$  being a cycle or a forest.
\end{abstract}

\section{Introduction}

All graphs in this paper are finite and simple. For a graph $G$, denote by $V(G)$, $E(G)$, $\Delta(G)$  and  $\bar d(G)=2|E(G)|/|V(G)|$  its {\it vertex set},   {\it edge set},  {\it maximum degree}, and {\it average degree}, respectively. Additionally, we say $G$ is of {\it order} $|V(G)|$.  For $v\in V(G)$, denote by  $d_G(v)$ and  $N_G(v)$ its {\it degree} and {\it neighbourhood}, respectively.

The study of graph inversions originated with the work of Belkhechine~\cite{belkhechine:thesis}. 
Given an oriented graph $\overrightarrow G$ and a set $X$ of its vertices, the {\it inversion} of $X$ reverses all arcs with both endpoints in $X$,    and denote the obtained oriented graph by  $\Inv(\overrightarrow G,X)$.
Much research has focused on the {\it inversion number} $\operatorname{inv}(D)$ of a digraph $D$, the minimum number of inversions required to render $D$ acyclic.
Belkhechine et al.~\cite{BBBP} proved that for every fixed $k$, deciding whether the inversion number of a tournament is at most $k$ is polynomial-time solvable; in contrast, Bang-Jensen et al.~\cite{BangJensenDaSilvaHavet} proved that deciding whether the inversion number of a digraph is $1$ is NP-complete.
For the maximum inversion number $\operatorname{inv}(n)$ over all oriented graphs of order $n$, Aubian et al.~\cite{AubianEtAl} and Alon et al.~\cite{AlonEtAl} independently proved 
 $n-2\sqrt{n\log n}\le \operatorname{inv}(n)\le n-\lceil\log_2(n+1)\rceil$. 
 Restricting the size of the sets, for a digraph $D$, denote by $\operatorname{inv}_k(D)$  the minimum number of inversions using sets of size  at most $k$
 that render $D$ acyclic. Let   $\operatorname{inv}_k(n)$ be the maximum inversion number  over all oriented graphs of order $n$.    
 Spencer~\cite{Spencer1971}  proved that $ \operatorname{inv}_2(n)\le n(n-1)/4-\Omega(n^{3/2})$. For more results about $ \operatorname{inv}_k(n)$, we refer to \cite{Chung1991,Spencer1980,Yuster2025,Vega1983}.   Beyond acyclicity,  Duron et al.~\cite{DuronEtAl} investigated the minimum number of inversions making a digraph $k$-arc-strong or $k$-strong.

The inversion diameter of a graph was introduced by Havet et al.~\cite{HavetHorschRambaud}.
Rather than reaching an acyclic orientation, it concerns transformations between arbitrary orientations of a fixed graph $G$.
The {\it inversion graph} $I(G)$ has as its vertices the labelled orientations of $G$; two orientations $\overrightarrow G_1$ and $\overrightarrow G_2$ are adjacent if there is a set $X$ of vertices with
$\Inv(\overrightarrow G_1,X)=\overrightarrow G_2$.
The {\it inversion diameter} $\diam(I(G))$ is the diameter of $I(G)$: the worst-case number of inversions needed to transform one orientation of $G$ into another.

Havet et al.~\cite{HavetHorschRambaud} showed that the inversion diameter is functionally equivalent to the star chromatic number $\chi_s(G)$, the acyclic chromatic number $\chi_a(G)$ and the oriented chromatic number $\chi_o(G)$, so that graph classes with any of these parameters bounded also have bounded inversion diameter.
They proved the following.

\begin{theorem}[{\cite{HavetHorschRambaud}}]\label{thm:HHR}
For every graph $G$, it holds  $\diam(I(G))\le 2\Delta(G)-1$.
 
\end{theorem}

They also showed $\diam(I(G))\le \Delta(G)$ when $\Delta(G)\le2$, and conjectured the following.

\begin{conjecture}[{\cite{HavetHorschRambaud}}]\label{conj:Delta}
For every graph $G$, it holds
 $\diam(I(G))\le \Delta(G)$.

\end{conjecture}

Wang et al.~\cite{WangWangYangLu} verified Conjecture~\ref{conj:Delta} for $\Delta(G)=3$ using computer assistance.
For bipartite graphs the stronger bound
$\diam(I(G))\le \Delta(G)+\lceil\log_2\Delta(G)\rceil-1$
was obtained in~\cite{HavetHorschRambaud}, and Arana et al.~\cite{AranaEtAl} extended such  results to all triangle-free graphs $G$: 
\begin{equation}\label{eq:triangle-free-delta-bound}
  \diam(I(G))\le \Delta(G)+\lfloor\log_2\Delta(G)\rfloor.
\end{equation}
It is natural to ask how the density of $G$ governs $\diam(I(G))$.
Havet et al.~\cite{HavetHorschRambaud} raised the following problem.

\begin{problem}[{\cite{HavetHorschRambaud}}]\label{prob:edge-bound}
Is there an absolute constant $\alpha$ such that
\[
  \diam(I(G))\le \alpha\sqrt{|E(G)|}
\]
for every graph $G$?
\end{problem}

This paper studies bounds dependent on edge number. 
Our first result (Theorem \ref{thm:edge-upper}) answers Problem~\ref{prob:edge-bound} affirmatively.

\begin{theorem}\label{thm:edge-upper}
For every graph $G$,
\[
  \diam(I(G))\le 2\sqrt{|E(G)|}.
\]
\end{theorem}

We show that the exponent $1/2$ is optimal in the strongest possible sense, via a clean lower bound that holds for every graph.
In Section~\ref{sec:preliminaries}, we   explains that, after fixing a reference orientation, $I(G)$ is a Cayley graph on $\F_2^{E(G)}$, which makes the bound intuitively clear.

\begin{theorem}\label{thm:lower-bound}
For every graph $G$ with at least one edge,
\[
  \diam(I(G))\ge \frac{|E(G)|}{|V(G)|}=\frac{\bar d(G)}2 .
\]
\end{theorem}

Combining Theorems~\ref{thm:edge-upper} and~\ref{thm:lower-bound} sandwiches the inversion diameter between $|E(G)|/|V(G)|$ and $2\sqrt{|E(G)|}$, both of which are $\Theta(n)$ for the complete graph $K_n$.
Thus $\sqrt{|E(G)|}$ is the correct order for every dense family, not merely for complete graphs.

We next characterize the optimal constant in Theorem~\ref{thm:edge-upper} which   drops to $\sqrt2$, once Conjecture~\ref{conj:Delta} holds.

\begin{theorem}\label{thm:optimal-constant}
Let $\alpha^*$ be the least constant such that
\[
  \diam(I(G))\le \alpha^*\sqrt{|E(G)|}
\]
for every graph $G$.
Then $\sqrt2\le \alpha^*\le2$, where $\sqrt2$ is witnessed by complete graphs,  with $\alpha^*=\sqrt2$ if Conjecture~\ref{conj:Delta} holds.
Moreover,
\[
  \diam(I(G))\le \sqrt{2|E(G)|}
\]
holds whenever $\diam(I(H))\le \Delta(H)$ for all subgraphs $H$ of $G$; in particular it holds for all $G$ with $\Delta(G)\le3$.
\end{theorem}

Furthermore, for triangle-free graphs, which include all bipartite graphs, the upper constant $\sqrt2$ is attained  up to a logarithmic term, by feeding~\eqref{eq:triangle-free-delta-bound} into the same induction.

\begin{theorem}\label{thm:triangle-free-edge}
For every triangle-free graph $G$ with $|E(G)|\ge1$,
\[
  \diam(I(G))\le \sqrt{2|E(G)|}+\frac12\log_2\bigl(2|E(G)|\bigr).
\]
\end{theorem}

We next refine the bipartite results.
Write $G[A,B]$ for a bipartite graph with parts $A$ and $B$, and let
\[
  \Delta(A):=\max_{a\in A}d_G(a),\qquad
  \Delta(B):=\max_{b\in B}d_G(b).
\]

\begin{theorem}\label{thm:bipartite-delta-three}
Let $G[A,B]$ be a bipartite graph with $\Delta(G)\le3$.
Then $\diam(I(G))\le \Delta(G)$, that is, Conjecture~\ref{conj:Delta} holds for bipartite graphs of maximum degree at most three.
\end{theorem}

\begin{theorem}\label{thm:asymmetric-bipartite}
Let $G[A,B]$ be a bipartite graph with $\Delta(A)=\sigma$ and $\Delta(B)=\rho$.
If $\rho=0$, then $\diam(I(G))=0$.
If $\rho\ge1$, then
\[
  \diam(I(G))\le
  \max\left\{\rho,
  \left\lceil\log_2\bigl(2+\sigma(2^{\rho-1}-1)\bigr)\right\rceil\right\}.
\]
The same bound holds with $\sigma$ and $\rho$ interchanged.
\end{theorem}

As a consequence, the inversion diameter of a bipartite graph can be exponentially smaller than its maximum degree (Corollary~\ref{cor:exp-smaller}).

Finally, we study lexicographic products with edgeless graphs.
For a graph $G$ and an integer $t\ge1$, denote by $G[\overline{K_t}]$  the {\it lexicographic product graph} obtained by replacing each vertex of $G$ with an independent set of size $t$ and joining two such sets completely whenever the corresponding vertices are adjacent in $G$. Here $\overline{K_t}$ is the edgeless graph on $t$ vertices.
Havet et al.~\cite{HavetHorschRambaud} conjectured that
$\diam(I(G[\overline{K_t}]))\le t\cdot\diam(I(G))$ (Conjecture~\ref{conj:lex-product}).
We confirm this for forests and cycles.

\begin{theorem}\label{thm:cycles}
For  every positive integer $t$, if $F$ is a forest or a cycle, then 
     $$
    \diam(I(F[\overline{K_t}]))\le t\cdot \diam(I(F)).
      $$
\end{theorem}

The paper is organised as follows.
In section~\ref{sec:preliminaries}, we develop the Cayley-graph viewpoint and the linear-algebraic reduction lemmas.
In section~\ref{sec:edge-bounds}, we prove bounds in terms of edge number for general graphs (Theorems~\ref{thm:edge-upper}, \ref{thm:lower-bound} and~\ref{thm:optimal-constant}).
In section~\ref{sec:triangle-free-bipartite}, we revise the bounds for triangle-free and bipartite graphs (Theorems~\ref{thm:triangle-free-edge}, \ref{thm:bipartite-delta-three} and~\ref{thm:asymmetric-bipartite}).
In section~\ref{sec:lex-products}, we deal with lexicographic products (Theorem~\ref{thm:cycles}).
%Section~\ref{sec:concluding} collects open problems.

\section{Preliminaries}\label{sec:preliminaries}

We first show that the inversion graph could be seen as a Cayley graph.  Fix a reference orientation $D_0$ of $G$ and identify each orientation $D$ with the vector
$\delta(D)\in\F_2^{E(G)}$ recording, for each edge, whether $D$ and $D_0$ disagree there, i.e.  for each edge $e$, the entry is 1 if they disagree on $e$, otherwise, 0.
For $X\subseteq V(G)$, let $\chi_X\in\F_2^{E(G)}$ be the indicator of the edges with both endpoints in $X$,
\[
  \chi_X(e)=1 \quad\Longleftrightarrow\quad \text{both endpoints of }e\text{ lie in }X .
\]
Inverting $X$ reverses exactly the edges counted by $\chi_X$, so
$\delta(\Inv(D,X))=\delta(D)+\chi_X$.

\begin{proposition}\label{prop:cayley}
Let
\[
  S=\{\chi_X:X\subseteq V(G)\}\subseteq\F_2^{E(G)}.
\]
If Cayley graphs are allowed to have loops, then
\[
  I(G)\cong \Cay(\F_2^{E(G)},S).
\]
Equivalently, the usual loopless inversion graph is obtained by using the connection set $S\setminus\{0\}$; deleting loops does not change distances.
Consequently $I(G)$ is vertex-transitive, and $\diam(I(G))$ equals the covering radius of $S$: the least $r$ such that every vector of $\F_2^{E(G)}$ is a sum of at most $r$ elements of $S$.
\end{proposition}

\begin{proof}
Under $\delta$, two orientations $D_1,D_2$ are adjacent in $I(G)$ if and only if
$\delta(D_2)-\delta(D_1)=\chi_X$ for some $X$, equivalently
$\delta(D_1)+\delta(D_2)\in S$, since we are over $\F_2$.
This is exactly the adjacency relation in the Cayley graph of the group $\F_2^{E(G)}$ with connection set $S$, apart from possible loops arising from $0\in S$.
The loopless version has the same distances.
A Cayley graph is vertex-transitive, and the distance from the identity to a vector $v$ is the least number of elements of $S$ summing to $v$.
Taking the maximum over $v$ gives the stated covering radius.
\end{proof}

The next two lemmas convert inversion diameter bounds into bilinear representation problems over $\F_2$. 
For $x,y\in\F_2^t$ write
$x\cdot y=\sum_{i=1}^t x_i y_i\in\F_2$
for the standard bilinear form.

\begin{lemma}\label{lem:bipartite-linear}
Let $G[A,B]$ be a bipartite graph and let $t\ge0$ be an integer.
Suppose that for every map $\pi:E(G)\to\F_2$,  there exist maps
$f:A\to\F_2^t$ and $g:B\to\F_2^t$ such that
\[
  f(a)\cdot g(b)=\pi(ab)
\]
for every $ab\in E(G)$.
Then $\diam(I(G))\le t$.
\end{lemma}

\begin{proof}
Let $\overrightarrow G_1,\overrightarrow G_2$ be arbitrary orientations and define $\pi(ab)=1$ if $\overrightarrow G_1$ and $\overrightarrow G_2$ orient $ab$ oppositely, and $0$ otherwise.
Choose $f,g$ as in the hypothesis and write
$f(a)=(f_1(a),\ldots,f_t(a))$ and $g(b)=(g_1(b),\ldots,g_t(b))$.
For $i\in\{1,\ldots,t\}$, set
\[
  X_i:=\{a\in A:f_i(a)=1\}\cup\{b\in B:g_i(b)=1\}.
\]
For an edge $ab$, the inversion on $X_i$ reverses $ab$ if and only if $f_i(a)=g_i(b)=1$.
Hence the parity of the number of times $ab$ is reversed by $X_1,\ldots,X_t$ is
\[
  \sum_{i=1}^t f_i(a)g_i(b)=f(a)\cdot g(b)=\pi(ab)\pmod2.
\]
Thus  the edges on which $\overrightarrow G_1$ and $\overrightarrow G_2$ disagree are reversed exactly an odd number of times, so applying the inversions $X_1,
\ldots,X_t$ to $\overrightarrow G_1$ yields $\overrightarrow G_2$.
As the two orientations were arbitrary, $\diam(I(G))\le t$.
\end{proof}

\begin{lemma}\label{lem:graph-linear}
Let $G$ be a graph and let $t\ge0$ be an integer.
Suppose that for every map $\pi:E(G)\to\F_2$, there is a map
$h:V(G)\to\F_2^t$ such that
\[
  h(u)\cdot h(v)=\pi(uv)
\]
for every $uv\in E(G)$.
Then $\diam(I(G))\le t$.
\end{lemma}

\begin{proof}
The proof is identical to Lemma~\ref{lem:bipartite-linear}, taking
$X_i=\{v\in V(G):h_i(v)=1\}$.
The parity of times that an edge $uv$ is reversed is 
$\sum_{i=1}^t h_i(u)h_i(v)=h(u)\cdot h(v)=\pi(uv)$ modulo $2$.
\end{proof}

The following monotonicity property of~\cite{HavetHorschRambaud} and  elementary inequality will also be applied later.

\begin{lemma}[{\cite[Lemma~2.1]{HavetHorschRambaud}}]\label{lem:vertex-deletion}
For every graph $G$ and every $v\in V(G)$,
\[
  \diam(I(G))\le \diam(I(G-v))+1.
\]
\end{lemma}

\begin{lemma}\label{lem:sqrt-inequality}
Let $m\ge1$ and $0\le\Delta\le m$.
If $\Delta>\sqrt{2m}$, then $\sqrt{2m}-\sqrt{2(m-\Delta)}>1$.
 
\end{lemma}

\begin{proof}
Using $\sqrt{2(m-\Delta)}\le\sqrt{2m}$, we have 
\[
\begin{aligned}
\sqrt{2m}-\sqrt{2(m-\Delta)}
  &=\frac{2\Delta}{\sqrt{2m}+\sqrt{2(m-\Delta)}}  \\
  &\ge \frac{2\Delta}{2\sqrt{2m}}
   =\frac{\Delta}{\sqrt{2m}}>1.
\end{aligned}
\]
\end{proof}

Applying Lemmas~\ref{lem:vertex-deletion} and \ref{lem:sqrt-inequality}, the following result holds.  
\begin{proposition}\label{prop:hereditary-sqrt2}
Let $\mathcal C$ be a class of graphs closed under vertex deletion such that every $H\in\mathcal C$ satisfies
$\diam(I(H))\le\Delta(H)$.
Then every $G\in\mathcal C$ satisfies
\[
  \diam(I(G))\le\sqrt{2|E(G)|}.
\]
\end{proposition}

\begin{proof}
We prove by inducting on $|V(G)|$. Set  $m=|E(G)|$ and $\Delta=\Delta(G)$.
If $m=0$, both sides are $0$.
If $\Delta\le\sqrt{2m}$, then
\[
  \diam(I(G))\le\Delta\le\sqrt{2m}
\]
since $G\in\mathcal C$.
Otherwise $\Delta>\sqrt{2m}$.
As $m\ge\Delta$, the graph $G'=G-v$ obtained by deleting a vertex of degree $\Delta$ has
$|E(G')|=m-\Delta\ge0$ and lies in $\mathcal C$.
By Lemma~\ref{lem:vertex-deletion}, the induction hypothesis and Lemma~\ref{lem:sqrt-inequality},
\[
  \diam(I(G))\le \diam(I(G'))+1
  \le \sqrt{2(m-\Delta)}+1<\sqrt{2m}.
\]
\end{proof}

\section{Bounds for general graphs}\label{sec:edge-bounds}
 
 In this section, we  prove bounds on  the inversion diameter for general graphs, specifically,  establish  Theorems~\ref{thm:edge-upper}, \ref{thm:lower-bound} and~\ref{thm:optimal-constant}. Let us  start from the upper bound. 

\begin{proof}[Proof of Theorem~\ref{thm:edge-upper}]
Write $m=|E(G)|$ and $\Delta=\Delta(G)$.
The claim is trivial when $m=0$, so assume $m\ge1$.
If $m\ge\Delta^2$, then $\Delta\le\sqrt m$ and Theorem~\ref{thm:HHR} gives
\[
  \diam(I(G))\le2\Delta-1\le2\sqrt m-1<2\sqrt m.
\]
Now suppose $m<\Delta^2$, i.e. $\Delta>\sqrt m$.
We induct on $|V(G)|$.
The base cases are immediate.
Let $v$ be a vertex of degree $\Delta$ and put $G'=G-v$.
Then $|E(G')|=m-\Delta$.
By the induction hypothesis,
\[
  \diam(I(G'))\le2\sqrt{m-\Delta}.
\]
By Lemma~\ref{lem:vertex-deletion},
\[
  \diam(I(G))\le \diam(I(G'))+1\le2\sqrt{m-\Delta}+1.
\]
Since $\Delta>\sqrt m$,
\[
  2\sqrt{m-\Delta}+1<2\sqrt{m-\sqrt m}+1.
\]
It remains only to observe that
\[
  2\sqrt{m-\sqrt m}+1<2\sqrt m,
\]
which is equivalent, after moving $1$ to the other side and squaring, to
$4m-4\sqrt m<4m-4\sqrt m+1$.
Thus $\diam(I(G))<2\sqrt m$, completing the induction.
\end{proof}

 Now, we move to the proof of the lower bound. 

\begin{proof}[Proof of Theorem~\ref{thm:lower-bound}]
Write $n=|V(G)|$, $m=|E(G)|$ and $k=\diam(I(G))$.
By Proposition~\ref{prop:cayley}, every vector of $\F_2^m$ is a sum of at most $k$ elements of
$S=\{\chi_X:X\subseteq V(G)\}$.
Since $0\in S$,  the reachable vectors form the image of the map
\[
  S^k\longrightarrow\F_2^m,
  \qquad (s_1,\ldots,s_k)\longmapsto s_1+\cdots+s_k,
\]
whose domain has $|S|^k\le(2^n)^k=2^{nk}$ elements.
As all $2^m$ vectors must be reached, $2^{nk}\ge2^m$, i.e. $nk\ge m$.
Hence
\[
  k\ge \frac mn=\frac{|E(G)|}{|V(G)|}=\frac{\bar d(G)}2 .
\]
\end{proof}

%Because $\diam(I(G))$ is an integer, Theorem~\ref{thm:lower-bound} gives
%$\diam(I(G))\ge\lceil |E(G)|/|V(G)|\rceil$.
%For $K_n$, this  yields  $\diam(I(K_n))\ge\lceil(n-1)/2\rceil$, within a factor $2$ of the true value $n-1$.

We finish this section by showing the bounds on the optimal constant, i.e.  Theorem~\ref{thm:optimal-constant}.   
\begin{proof}[Proof of Theorem~\ref{thm:optimal-constant}]
The bound $\alpha^*\le2$ is Theorem~\ref{thm:edge-upper}.
For  result under the assumption that Conjecture~\ref{conj:Delta} holds and the last sentence, apply Proposition~\ref{prop:hereditary-sqrt2}: the class of all graphs satisfies its hypothesis under Conjecture~\ref{conj:Delta}, giving $\alpha^*=\sqrt2$; the class $\{G:\Delta(G)\le3\}$ is closed under vertex deletion and satisfies the hypothesis, by~\cite{HavetHorschRambaud} for $\Delta\le2$ and by~\cite{WangWangYangLu} for $\Delta=3$, giving $\diam(I(G))\le\sqrt{2|E(G)|}$ for $\Delta(G)\le3$.
Finally, for the lower bound, the complete graph $K_n$ satisfies
\[
  \frac{\diam(I(K_n))}{\sqrt{|E(K_n)|}}
  =\frac{n-1}{\sqrt{\binom n2}}
  =\sqrt{\frac{2(n-1)}n}\longrightarrow\sqrt2,
\]
where $\diam(I(K_n))=n-1$ is from~\cite{HavetHorschRambaud}  and so $\alpha^*\ge\sqrt2$.
\end{proof}

\section{Triangle-free and bipartite graphs}\label{sec:triangle-free-bipartite}

In this section, we revise the upper bound   for triangle-free graphs, especially bipartite graphs. For triangle-free graphs, the base case~\eqref{eq:triangle-free-delta-bound} is within a logarithmic term of Conjecture~\ref{conj:Delta}.

\begin{proof}[Proof of Theorem~\ref{thm:triangle-free-edge}]
Set
\[
  \varphi(m)=\sqrt{2m}+\frac12\log_2(2m)
\]
and let $m=|E(G)|\ge1$ and $\Delta=\Delta(G)\ge1$. We prove by inducting  on $|V(G)|$. 
If $\Delta\le\sqrt{2m}$, then by~\eqref{eq:triangle-free-delta-bound} and monotonicity of $\log_2$,
\[
\begin{aligned}
  \diam(I(G))
  &\le \Delta+\lfloor\log_2\Delta\rfloor \\
  &\le \sqrt{2m}+\log_2\sqrt{2m}
   = \sqrt{2m}+\frac12\log_2(2m)=\varphi(m).
\end{aligned}
\]
If $\Delta>\sqrt{2m}$, by deleting a vertex of degree $\Delta$, we obtain a triangle-free graph $G'$ with
$m' =|E(G')|=m-\Delta\ge0$.
If $m'=0$, then $\diam(I(G))\le1\le\varphi(m)$.
If $m'\ge1$, then by Lemma~\ref{lem:vertex-deletion}, the induction hypothesis and Lemma~\ref{lem:sqrt-inequality},
\[
\begin{aligned}
  \varphi(m)-\varphi(m')
  &=\bigl(\sqrt{2m}-\sqrt{2m'}\bigr)
    +\frac12\bigl(\log_2(2m)-\log_2(2m')\bigr) \\
  &\ge \sqrt{2m}-\sqrt{2m'}>1,
\end{aligned}
\]
since the logarithmic difference is nonnegative.
Hence
\[
  \diam(I(G))\le\varphi(m')+1<\varphi(m) 
\]
and so we are done. 
\end{proof}

\begin{remark}\label{rem:triangle-free-constant}
The constant $\sqrt2$ cannot be obtained by this method with $\sqrt{2m}$ replaced by $\sqrt m$: when $\Delta\approx\sqrt m$, deleting one maximum-degree vertex decreases $\sqrt m$ by only about $1/2$.
Determining the optimal triangle-free constant remains open, but balanced complete bipartite graphs give the lower bound $1$.
Indeed, for $K_{a,b}$, after fixing a reference orientation, orientation differences are $a\times b$ matrices over $\F_2$.
An inversion toggles a matrix of the form $xy^T$, where $x\in\F_2^a$ and $y\in\F_2^b$, hence a matrix of rank at most $1$.
Conversely every matrix of rank $r$ over $\F_2$ is a sum of $r$ rank-one matrices.
Therefore the distance from $0$ to a matrix $M$ is $\operatorname{rank}(M)$, and so
\[
  \diam(I(K_{a,b}))=\min\{a,b\}.
\]
Taking $a=b=r$ gives $|E(K_{r,r})|=r^2$ and $\diam(I(K_{r,r}))=r=\sqrt{|E(K_{r,r})|}$.
Thus the least asymptotic triangle-free constant $\beta$ satisfies  $1\le\beta\le\sqrt2$ .
 
\end{remark}

The proof of Theorem~\ref{thm:bipartite-delta-three} uses Lemma~\ref{lem:bipartite-linear} together with the following representation theorem.

\begin{theorem}[\cite{FurediEtAl}] \label{thm:FGHK}
Let $H$ be a $3$-uniform hypergraph with maximum degree at most $3$.
Then there is a map $c:V(H)\to\F_2^3$ such that for every hyperedge
$\{x,y,z\}\in E(H)$, the vectors $c(x),c(y),c(z)$ form a basis of $\F_2^3$.
\end{theorem}

\begin{proof}[Proof of Theorem~\ref{thm:bipartite-delta-three}]
 If $\Delta(G)\le 2$, then by Corollary 4.4 for forests and Theorem 6.4 for graphs $G$ with $\Delta(G)= 2$ in \cite{HavetHorschRambaud}, we are done. Thus suppose $\Delta(G)=3$.
Build a $3$-uniform hypergraph $H$ on $A$ together with auxiliary vertices as follows: for each $b\in B$ with $d_G(b)\ge1$, set
\[
  E_b:=
  \begin{cases}
    N_G(b), & d_G(b)=3,\\
    N_G(b)\cup\{z_b\}, & d_G(b)=2,\\
    N_G(b)\cup\{z_b,z'_b\}, & d_G(b)=1,
  \end{cases}
\]
where the auxiliary vertices are new.
Let $E(H)$ be the distinct $3$-sets among the $E_b$.
Then each $a\in A$ lies in at most $d_G(a)\le3$ hyperedges and each auxiliary vertex lies in exactly one, so $\Delta(H)\le3$.
By Theorem~\ref{thm:FGHK}, there is $c:V(H)\to\F_2^3$ making the vectors on each hyperedge a basis. Set $f(a)=c(a)$ for each $a\in A$. 
Then $\{f(a):a\in N_G(b)\}$ is linearly independent, as a subset of the basis on $E_b$, for every $b$ with $d_G(b)\ge1$. For each $b$ with $d_G(b)\ge1$, the map
\[
  T_b:\F_2^2\to\F_2^{N_G(b)},\qquad
  T_b(y)=(f(a)\cdot y)_{a\in N_G(b)},
\]
is surjective, yielding $g(b)\in\F_2^3$ with $f(a)\cdot g(b)=\pi(ab)$. Set $g(b)=0$ if $d_G(b)=0$. Then $\diam(I(G))\le3$ by
Lemma~\ref{lem:bipartite-linear}.
\end{proof}

Finally, we provide the upper bound on inversion diameter of bipartite graphs in terms of vertex degree. 

\begin{proof}[Proof of Theorem~\ref{thm:asymmetric-bipartite}]
If $\rho=0$, then $E(G)=\emptyset$ and $\diam(I(G))=0$.
Assume $\rho\ge1$ and let
\[
  t=\max\left\{\rho,
  \left\lceil\log_2\bigl(2+\sigma(2^{\rho-1}-1)\bigr)\right\rceil\right\},
\]
so that $t\ge\rho$ and
\[
  2^t>1+\sigma(2^{\rho-1}-1).
\]
By Lemma~\ref{lem:bipartite-linear} it suffices to represent every $\pi:E(G)\to\F_2$ as
$\pi(ab)=f(a)\cdot g(b)$ with $f:A\to\F_2^t$ and $g:B\to\F_2^t$. Now, we proceed by defining $f$ and $g$. 

\smallskip
\noindent\emph{Choosing $f$, independently of $\pi$.}
Write $A=\{a_1,\ldots,a_n\}$.
We choose $f(a_1),f(a_2),\ldots$ greedily, maintaining the invariant that for every $b\in B$ the already chosen vectors among $\{f(a):a\in N_G(b)\}$ are linearly independent.
Take $f(a_1)$ to be any nonzero vector.
Suppose $f$ is defined on $a_1,\ldots,a_{i-1}$.
For each $b\in N_G(a_i)$, let
\[
  W_b=\spanop\{f(a_j):a_j\in N_G(b),\ j<i\}.
\]
Choosing $f(a_i)\notin\bigcup_{b\in N_G(a_i)} W_b$ preserves the invariant.
At most $d_G(b)-1\le\rho-1$ vectors of $N_G(b)$ are already chosen, and they are independent by the invariant, so
$\dim W_b\le\rho-1$ and $|W_b|\le2^{\rho-1}$.
Since each $W_b$ contains $0$ and $|N_G(a_i)|\le\sigma$,
\[
  \left|\bigcup_{b\in N_G(a_i)} W_b\right|
  \le 1+\sum_{b\in N_G(a_i)}(|W_b|-1)
  \le1+\sigma(2^{\rho-1}-1)<2^t=|\F_2^t|,
\]
so a valid $f(a_i)$ exists.
Proceeding to the end yields $f$ with $\{f(a):a\in N_G(b)\}$ independent for every $b$.

\smallskip
\noindent\emph{Solving for $g$.}
Fix $\pi$.
For each $b$ with $d_G(b)\ge1$ the map
\[
  T_b:\F_2^t\to\F_2^{N_G(b)},\qquad
  T_b(y)=(f(a)\cdot y)_{a\in N_G(b)},
\]
has independent rows, using $t\ge\rho\ge|N_G(b)|$, hence is surjective.
Pick $g(b)$ with $f(a)\cdot g(b)=\pi(ab)$ for all $a\in N_G(b)$.
Set $g(b)=0$ if $d_G(b)=0$.
By Lemma~\ref{lem:bipartite-linear}, $\diam(I(G))\le t$.
Interchanging $A$ and $B$ gives the symmetric bound.
\end{proof}

Taking $\rho=3$, we have 
$\lceil\log_2(2+3\sigma)\rceil\le\sigma$ for all $\sigma\ge4$. Combined with Theorem~\ref{thm:bipartite-delta-three}, which covers $\sigma\le3$, and Theorem~\ref{thm:asymmetric-bipartite}, we obtain the following extension of Theorem~\ref{thm:bipartite-delta-three}.

\begin{corollary}\label{cor:one-side-three}
Let $G[A,B]$ be a bipartite graph with $\Delta(B)=3$ and $\Delta(A)$ arbitrary.
Then $\diam(I(G))\le\Delta(G)$.
\end{corollary}

Generally when the two parts have very different maximum degrees, the result in Theorem~\ref{thm:asymmetric-bipartite} is far stronger than the maximum degree bounds.  

\begin{corollary}\label{cor:exp-smaller}
If one part of a bipartite graph $G$ has bounded degree $\rho=O(1)$, then
\[
  \diam(I(G))=O(\log\Delta(G)).
\]
%In particular, for $K_{2,N}$ with $N\ge2$,
%\[
  %\diam\bigl(I(K_{2,N})\bigr)\le\lceil\log_2(N+2)\rceil,
%\]
%which is exponentially smaller than $\Delta(K_{2,N})=N$.
\end{corollary}

\begin{proof}
Take the bounded side to be $B$, so $\rho=O(1)$ and $\sigma\le\Delta(G)$.
By Theorem~\ref{thm:asymmetric-bipartite}, we have 
\[
  \diam(I(G))\le
  \max\left\{\rho,
  \left\lceil\log_2\bigl(2+\sigma(2^{\rho-1}-1)\bigr)\right\rceil\right\}
  =O(\log_2\sigma)=O(\log_2\Delta(G)).
\]\end{proof}
%For $K_{2,N}$ we have $\sigma=N$ and $\rho=2$, so the bound is
%$\max\{2,\lceil\log_2(N+2)\rceil\}=\lceil\log_2(N+2)\rceil$ for $N\ge2$.

%By Theorem~\ref{thm:lower-bound},
%$\diam(I(K_{2,N}))\ge 2N/(N+2)\to2$, so the inversion diameter of $K_{2,N}$ lies between a constant and $\log_2(N+2)$. 

\section{Lexicographic products}\label{sec:lex-products}

For lexicographic product,  
Havet et al.~\cite{HavetHorschRambaud} proposed the following conjecture.

\begin{conjecture}[{\cite{HavetHorschRambaud}}]\label{conj:lex-product}
For every graph $G$ and every positive integer $t$,
\[
  \diam(I(G[\overline{K_t}]))\le t\cdot\diam(I(G)).
\]
\end{conjecture}

To attack such products, we apply a matrix version of Lemma~\ref{lem:graph-linear}.   A $t\times t$ matrices over $\F_2$ will be denoted by $\F_2^{t\times t}$.

\begin{lemma}\label{lem:matrix-product}
Let $G$ be a graph and let $t,\ell$ be positive integers.
Fix an arbitrary orientation of the edges of $G$.
Suppose that for every map
$\pi:E(G)\to\F_2^{t\times t}$ there exists a map
$U:V(G)\to\F_2^{\ell\times t}$ such that
\[
  \pi(uv)=U(u)^T U(v)
\]
for every oriented edge $uv$ of $G$.
Then
\[
  \diam(I(G[\overline{K_t}]))\le \ell.
\]
\end{lemma}

\begin{proof}
Denote the associated vertices in $G[\overline{K_t}]$ of $v\in V(G)$ by  $\{(v,1),\ldots,(v,t)\}$.
For two orientations of $G[\overline{K_t}]$, define for each oriented edge $uv$ of $G$ a matrix
$\pi(uv)\in\F_2^{t\times t}$ whose $(i,j)$ entry is $1$ exactly when the two orientations disagree on the edge $(u,i)(v,j)$.
Choose $U$ as in the hypothesis.
For $s\in\{1,\ldots,\ell\}$, set
\[
  X_s:=\{(v,i): U(v)_{s,i}=1\}.
\]
For an edge $(u,i)(v,j)$ lying above an oriented base edge $uv$, the inversion on $X_s$ reverses it if and only if
$U(u)_{s,i}=U(v)_{s,j}=1$.
Thus the parity of the number of reversals over $s=1,
\ldots,\ell$ is
\[
  \sum_{s=1}^{\ell} U(u)_{s,i}U(v)_{s,j}
  =\bigl(U(u)^T U(v)\bigr)_{i,j}
  =\pi(uv)_{i,j}.
\]
Hence the $\ell$ inversions transform one orientation into the other.
Since the orientations were arbitrary, the diameter is at most $\ell$.
\end{proof}

\begin{lemma}\label{lem:two-equations}
Let $\K$ be a field and let $B,C\in \K^{2t\times t}$.
If the concatenation $(B\mid C)\in \K^{2t\times2t}$ is invertible, then for all $X,Y\in \K^{t\times t}$ there is a unique $Z\in \K^{2t\times t}$ with
\[
  B^T Z=X
  \qquad\text{and}\qquad
  Z^T C=Y.
\]
\end{lemma}

\begin{proof}
The second equation is equivalent to $C^T Z=Y^T$, so the pair is equivalent to
\[
  \begin{pmatrix}B^T\\ C^T\end{pmatrix}Z
  =
  \begin{pmatrix}X\\ Y^T\end{pmatrix}.
\]
The coefficient matrix equals $(B\mid C)^T$, which is invertible since $(B\mid C)$ is. Hence $Z$ exists and is unique.
\end{proof}

\begin{lemma}\label{lem:mat}
    Let $\K$ be a field and let $M_1,M_2,...,M_k\in \K^{t\times t}$ with $k\ge 3$. Then there exist  $A_1, A_2,...,A_k\in \K^{2t\times t}$ such that 
    $$
    A_1^TA_2=M_1,~ A_2^TA_3=M_2,...,\text{ and } A_k^TA_1=M_k.
    $$
\end{lemma}

\begin{proof}

It is sufficient to consider the following two cases. 
\medskip

\noindent {\bf Case 1. } $k=2r+1$ where $r\ge1$.  

Put
\[
  P=\begin{pmatrix}I_t\\0\end{pmatrix},
  \qquad
  Q=\begin{pmatrix}M_k^T\\I_t\end{pmatrix},
  \qquad
  R=P+Q.
\]
Then $(P\mid Q)$ and $(P\mid R)$  are   invertible, since they are block-triangular with identity diagonal blocks. Moreover,   $(Q\mid R)$ is invertible. If $Qx+Ry=0$, then  since $R=P+Q$, one has
$Q(x+y)+Py=0$, hence $x+y=0$ and $y=0$ because $(P\mid Q)$ is invertible, so $x=y=0$.

Let $B_1=P$, let $B_{r+1}=Q$, and alternate $P,R$ for $B_1,
\ldots,B_r$, so $B_1=P,B_2=R,B_3=P,\ldots$.
Adjacent pairs $(B_i,B_{i+1})$ are always one of
$(P,R),(R,P),(P,Q),(R,Q)$, hence have invertible concatenation.
Let 
\[
  A_{2i-1}=B_i\qquad (i=1,\ldots,r+1),
\]
noting that $2(r+1)-1=k$.
For each $i=1,\ldots,r$, by Lemma~\ref{lem:two-equations}, there is  $A_{2i}$ with
\[
  B_i^T A_{2i}=M_{2i-1}
  \qquad\text{and}\qquad
  A_{2i}^T B_{i+1}=M_{2i},
\]
so that $A_j^T A_{j+1}=M_j$ for $j=1,\ldots,2r$.
Finally
\[
  A_{2r+1}^T A_1=Q^TP=M_k=M_{2r+1},
\]
since $Q^TP=(M_k\mid I_t)\binom{I_t}{0}=M_k$.

\medskip

\noindent {\bf Case 2. } $k=2r$ with $r\ge2$.

Put
\[
  P=\begin{pmatrix}I_t\\0\end{pmatrix},
  \qquad
  Q=\begin{pmatrix}0\\I_t\end{pmatrix},
  \qquad
  R=\begin{pmatrix}I_t\\I_t\end{pmatrix},
\]
for which $(P\mid Q)$, $(P\mid R)$ and $(Q\mid R)$ are clearly invertible.

Let $B_1=B_{r+1}=P$.
If $r$ is even, alternate $P,Q$ for $B_1,
\ldots,B_r$, so $B_r=Q$.
If $r$ is odd, alternate $P,Q$ for $B_1,
\ldots,B_{r-1}$ and set $B_r=R$.
In either case adjacent pairs $(B_i,B_{i+1})$ have invertible concatenation. 
Let 
\[
  A_{2i-1}=B_i\qquad (i=1,\ldots,r),
\]
and for each $i=1,\ldots,r$, by  Lemma~\ref{lem:two-equations}, there is  $A_{2i}$  with
\[
  B_i^T A_{2i}=M_{2i-1}
  \qquad\text{and}\qquad
  A_{2i}^T B_{i+1}=M_{2i}.
\]
Then $A_j^T A_{j+1}=M_j$ for $j=1,\ldots,2r-1$, and since
$B_{r+1}=B_1=A_1$ the wrap-around edge gives
$A_{2r}^T A_1=M_{2r}$.
\end{proof}

Now, we show that Conjecture~\ref{conj:lex-product} holds for forests and cycles. 
\begin{proof}[Proof of Theorem~\ref{thm:cycles}]

Suppose that $F=C_k$. Let  $C_k=v_1v_2\cdots v_kv_1$ with edges $e_i=v_iv_{i+1}$, indices modulo $k$.  Then $\diam(C_k)\ge 2$  by considering the distance between $v_1\to \cdots \to v_k\to v_1$ and   $v_1\to v_2\leftarrow v_3\leftarrow \cdots \leftarrow v_k\leftarrow v_1$.   
Let $\pi:E(C_k)\to\F_2^{t\times t}$ with $\pi(e_i)=M_i$.
By Lemma~\ref{lem:matrix-product} with $\ell=2t$ and Lemma~\ref{lem:mat},  there is
$f:V(C_k)\to\F_2^{2t\times t}$ with
\[
  f(v_j)^T f(v_{j+1})=M_j
\]
for all $j$ modulo $k$.   Suppose that $F$ is a forest, by Corollary 4.1 in \cite{HavetHorschRambaud}, the same argument works for $F$ being a forest with $\diam(I(F))=2$. If $E(F)=\emptyset$, then we are done. It remains to  assume $F=K_{1,t}$ whose $\diam(I(F))=1$. Let $v$ be the vertex that forms a single partite. By recursively considering oriented edges incident with $(v,i)$, we know $\diam(I(F[\overline{K_t}]))\le t$. 
\end{proof}

\begin{remark}\label{rem:cycle-lower}
The graph $C_k[\overline{K_t}]$ is $2t$-regular.
Moreover, $|V(C_k[\overline{K_t}])|=kt$ and $|E(C_k[\overline{K_t}])|=kt^2$, so Theorem~\ref{thm:lower-bound} gives
\[
  \diam(I(C_k[\overline{K_t}]))\ge t.
\]
Hence
\[
  t\le \diam\bigl(I(C_k[\overline{K_t}])\bigr)\le2t,
\]
and the bound of Theorem~\ref{thm:cycles} is optimal up to a factor of two.
\end{remark}

\section*{Acknowledgments}

  Jiawen Bo is supported by Fundamental and Interdisciplinary Disciplines Breakthrough Plan of the Ministry of Education of China (JYB2025XDXM207).
% \section{Concluding remarks}\label{sec:concluding}

% Recall that for every graph $G$ we have established
% \[
%   \frac{|E(G)|}{|V(G)|}\le \diam(I(G))\le2\sqrt{|E(G)|}. 
% \]
% Furthermore, the constant in the upper bound cannot be below $\sqrt2$, because of complete graphs. Several problems remain.

% \begin{problem}\label{ques:const}
% Is
% \[
%   \diam(I(G))\le\sqrt{2|E(G)|}
% \]
% for every graph $G$?
% \end{problem}
% By Proposition~\ref{prop:hereditary-sqrt2} this follows from Conjecture~\ref{conj:Delta}.  
% For triangle-free graphs and  lexicographic products, one may consider the followings. 

% \begin{problem}\label{prob:trianglefree}
%     What is the least $\beta$ such that 
% \[
%   \diam(I(G))\le(\beta+o(1))\sqrt{|E(G)|}
% \]
% for all triangle-free $G$?
% \end{problem}
% Recall that Theorem~\ref{thm:triangle-free-edge} and Remark~\ref{rem:triangle-free-constant} give $  1\le\beta\le\sqrt2$.   

% \begin{problem}\label{prob:Lexi}
%     Beyond  forests and  cycles  (Theorem \ref{thm:cycles}), does Conjecture~\ref{conj:lex-product} hold for all graphs, or at least with the constant in Theorem~\ref{thm:cycles} reduced from $2t$ towards the lower bound $t$ of Remark~\ref{rem:cycle-lower}?
% \end{problem}


\begin{thebibliography}{9}

\bibitem{AlonEtAl}
N.~Alon, E.~Powierski, M.~Savery, A.~Scott, and E.~Wilmer,
Invertibility of digraphs and tournaments,
\emph{SIAM J. Discrete Math.} 38(1) (2024), 327--347.

\bibitem{AranaEtAl}
C.~Arana, T.~Bellitto, H.~Buffi\`ere, Q.~Chuet, T.~Pierron, and A.~Reinald,
Inversion diameter and 2-edge-colored homomorphisms, \emph{arXiv:2602.24171} 
(2026).

\bibitem{AubianEtAl}
G.~Aubian, F.~Havet, F.~H\"orsch, F.~Klingelhoefer, N.~Nisse, C.~Rambaud, and Q.~Vermande,
Problems, proofs, and disproofs on the inversion number,
\emph{Electron. J. Combin.} 32(1) (2025), P1.42. 

\bibitem{BangJensenDaSilvaHavet}
J.~Bang-Jensen, J.~C.~F. da Silva, and F.~Havet,
On the inversion number of oriented graphs,
\emph{Discrete Math. Theor. Comput. Sci.} 23(2) (2022).


\bibitem{belkhechine:thesis}
H.~Belkhechine, 
   {Ind{\'e}composabilit{\'e} des graphes et des tournois}, \emph{
  Theses, {Universit{\'e} Claude Bernard-Lyon I; Universit{\'e} de
  Sfax. Facult{\'e} des sciences}},  (2009).

  
\bibitem{BBBP}
H.~Belkhechine, M.~Bouaziz, I.~Boudabbous, and M.~Pouzet,
Inversion dans les tournois,
\emph{C. R. Math. Acad. Sci. Paris} 348(13--14) (2010), 703--707.


\bibitem{Chung1991}
    F. R. Chung and R. L. Graham,  Quasi‐random tournaments, \emph{J. Graph Theory} 15(2) (1991), 173--198.

    
\bibitem{DuronEtAl}
J.~Duron, F.~Havet, F.~H\"orsch, and C.~Rambaud,
On the minimum number of inversions to make a digraph $k$-(arc-)strong,
\emph{J. Graph Theory} 111(2) (2025), 31--62.

\bibitem{FurediEtAl}
Z.~F\"uredi, J.~R. Griggs, R.~Holzman, and D.~J. Kleitman,
Representations of families of triples over GF(2),
\emph{J. Combin. Theory Ser. A} 53(2) (1990), 306--315.

%\bibitem{HavetHorschRambaud}
%F.~Havet, F.~H\"orsch, and C.~Rambaud,
%Diameter of the inversion graph, \emph{ arXiv:2405.04119v2} 
%(2024).

\bibitem{HavetHorschRambaud}
F.~Havet, F.~H\"orsch, and C.~Rambaud,
Diameter of the inversion graph,
\emph{Innov. Graph Theory} 3 (2026), 49--88.

\bibitem{Spencer1971}
J. Spencer,  Optimal ranking of tournaments,  \emph{Networks} 1(2) (1971),  135--138.

\bibitem{Spencer1980}
 J. Spencer,  Optimally ranking unrankable tournaments,  \emph{Period Math Hung} 11(2) (1980), 131--144.
 
\bibitem{Yuster2025}
 R. Yuster,  On tournament inversion,  \emph{J. Graph Theory}  110 (2025), 82--91.


 


\bibitem{Vega1983}
  W. F. de la Vega,  On the maximum cardinality of a consistent set of arcs in a random tournament,  \emph{J.  Combin. Theory Ser. B} 35(3) (1983),  328--332.
 
%\bibitem{WangWangYangLu}
%Y.~Wang, H.~Wang, Y.~Yang, and M.~Lu,
%Inversion diameter and treewidth, \emph{arXiv:2407.15384v4} 
%(2024).

\bibitem{WangWangYangLu}
Y.~Wang, H.~Wang, Y.~Yang, and M.~Lu,
Inversion diameter and treewidth,
\emph{Discrete Math. Theor. Comput. Sci.} 28(2) (2026).

\end{thebibliography}
\end{document}